\documentclass[12pt]{article}
\usepackage{amsmath}
\usepackage{amssymb}
\usepackage{amsthm}
\usepackage[mathscr]{eucal}
\usepackage{latexsym}
\usepackage{amscd}
\usepackage{epic}

\begin{document}
\baselineskip=15pt

\voffset -1.5truecm
\oddsidemargin .5truecm
\evensidemargin .5truecm

\theoremstyle{plain}
\swapnumbers\newtheorem{lema}{Lemma}[section]
\newtheorem{prop}[lema]{Proposition}
\newtheorem{coro}[lema]{Corolary}
\newtheorem{teor}[lema]{Theorem}
\newtheorem{demo}[lema]{Proof of Theorem 1}
\newtheorem{ejem}[lema]{Example}

\renewcommand{\refname}{\large\bf References}
\renewcommand{\thefootnote}{\fnsymbol{footnote}}

\def\flecha{\longrightarrow}
\def\asocia{\longmapsto}
\def\sii{\Longleftrightarrow}
\def\vector#1#2{({#1}_1,\dots,{#1}_{#2})}
\def\conjunto#1{\boldsymbol{[}\,{#1}\,\boldsymbol{]}}
\def\particion{\vdash}

\def\natural{{\mathbb N}}
\def\entero{{\mathbb Z}}
\def\racional{{\mathbb Q}}
\def\real{{\mathbb R}}
\def\complejo{{\mathbb C}}

\def\liri{Littlewood-Richardson\ }
\def\beze{Berenstein-Zelevinsky\ }
\def\cont{content\ }
\def\lr#1{{\sf LR}_{#1}}
\def\lrtipo#1{{\sf LR}_{#1}(\lambda,\mu,\nu)}
\def\pesos#1{{\sf D}_{#1}}
\def\triangulos#1{\real^{\binom{k+2}{2}-1}}
\def\tipo{(\lambda,\mu,\nu)}
\def\size#1{\vert {#1} \vert}
\def\arreglo#1{({#1}_{ij})_{\,0\le i\le j\le k}}
\def\hive#1{{\sf H}_{#1}}
\def\hivetipo#1{{\sf H}_{#1}(\lambda,\mu,\nu)}
\def\bz#1{{\sf BZ}_{#1}}
\def\bztipo#1{{\sf BZ}_{#1}(\lambda,\mu,\nu)}
\def\lrcoef{c^\lambda_{\mu\,\nu}}
\def\ent#1#2#3{{#1}_{{#2}\,{#3}}}

\def\carre#1{{\sf H}_{#1}(\lambda,\nu,\mu)}

\def\forma#1{{\sf sh}({#1})}
\def\sst{semistandard\ tableau\ }
\def\sucesion#1#2{{#1}_1 < {#1}_2 < \cdots < {#1}_{#2}}
\def\parti#1{\pi(#1)}

\def\HG{\Delta}
\def\BZ{\Gamma}
\def\rr{\mathbb{R}}
\def\vp{\varphi}


\setlength\unitlength{0.08em}
\savebox0{\rule[-2\unitlength]{0pt}{10\unitlength}%
\begin{picture}(10,10)(0,2)
\put(0,0){\line(0,1){10}}
\put(0,10){\line(1,0){10}}
\put(10,0){\line(0,1){10}}
\put(0,0){\line(1,0){10}}
\end{picture}}

\newlength\cellsize \setlength\cellsize{18\unitlength}
\savebox2{%
\begin{picture}(18,18)
\put(0,0){\line(1,0){18}}
\put(0,0){\line(0,1){18}}
\put(18,0){\line(0,1){18}}
\put(0,18){\line(1,0){18}}
\end{picture}}
\newcommand\cellify[1]{\def\thearg{#1}\def\nothing{}%
\ifx\thearg\nothing
\vrule width0pt height\cellsize depth0pt\else
\hbox to 0pt{\usebox2\hss}\fi%
\vbox to 18\unitlength{
\vss
\hbox to 18\unitlength{\hss$#1$\hss}
\vss}}
\newcommand\tableau[1]{\vtop{\let\\=\cr
\setlength\baselineskip{-16000pt}
\setlength\lineskiplimit{16000pt}
\setlength\lineskip{0pt}
\halign{&\cellify{##}\cr#1\crcr}}}
\savebox3{%
\begin{picture}(15,15)
\put(0,0){\line(1,0){15}}
\put(0,0){\line(0,1){15}}
\put(15,0){\line(0,1){15}}
\put(0,15){\line(1,0){15}}
\end{picture}}
\newcommand\expath[1]{%
\hbox to 0pt{\usebox3\hss}%
\vbox to 15\unitlength{
\vss
\hbox to 15\unitlength{\hss$#1$\hss}
\vss}}


\begin{centering}
{\Large\bf  Combinatorics and geometry of}\\[.2cm]
{\Large\bf  Littlewood-Richardson cones}\\[1cm]
{{\large\sf  }}\vskip.1cm
{{\large{\bf Igor Pak}
\vskip.1cm
Massachusetts Institute of Technology\\
Cambridge, MA 02138 USA\\
e-mail: pak@math.mit.edu\\}
\vskip.4cm
{\large\sf and}\vskip.4cm
{{\large{\bf Ernesto Vallejo} \footnotemark
\vskip.1cm
Instituto de Matem\'aticas, Unidad Morelia, UNAM\\
Apartado Postal 61-3, Xangari\\
58089 Morelia, Mich., MEXICO\\
e-mail: vallejo@matmor.unam.mx}}
\vskip.3cm
March 21, 2003}
\\
\end{centering}
\footnotetext{This work was done during a sabbatical stay at MIT
Mathematics Department. I would like to thank CONACYT and DGAPA-UNAM
for financial support.}

\vskip 2pc
\begin{abstract}
We present several direct bijections between different
combinatorial interpretations of the Littlewood-Richardson
coefficients.
The bijections are defined by explicit
linear maps which have other applications.
\end{abstract}

\vskip 2.5pc
{\large\bf 1\quad Introduction}

\vskip 1.5pc

In the past decade the Littlewood-Richardson rule (LR~rule)
moved into a center stage in the combinatorics of Young tableaux.
Much attention have received classical applications
(to representation theory of the symmetric and the full
linear group, to the symmetric functions, etc.) as well as more
recent developments (Schubert calculus, eigenvalues of
Hermitian matrices, etc.)  While various combinatorial
interpretations of the Littlewood-Richardson coefficients
were discovered, there seems to be little understanding
of how they are related to each other, and little order among
them. This paper makes a new step in this direction.

We start with three major combinatorial interpretations of the
LR coefficients which we view as integer points in certain
cones.  We present simple linear maps between the cones
which produce explicit bijections for all triples of partitions
involved in the LR rule.  These bijections are quite
natural in this setting and in a certain sense can be
shown to be unique. Below we further emphasize the
importance of the linear maps.

\smallskip

A classical version of the LR~rule, in terms
of certain Young tableaux, is now well understood, and
its proof has been perfected for decades.
We refer to~\cite{lee} for a beautifully written survey
of the ``classical'' approach, with a historical overview
and connections to the jeu-de-taquin, Sch\" utzenberger
involution, etc.  Unfortunately, the language of Young tableaux
is often too rigid to be able to demonstrate the inherent
symmetries of the LR~coefficients.

A radically different combinatorial interpretation in due
to Berenstein and Zelevinsky, in terms of the so called
BZ triangles, which makes explicit all but
one symmetry of the LR~coefficients\footnote{We should
warn the reader that the BZ triangles presented
in~\cite{stan} are different, but strongly related.}.
The authors' proof
in~\cite{beze} relies on a series of previous
papers~\cite{gz,bz1,bz2}, a situation that is
hardly satisfactory.  A paper~\cite{car} establishes a
technically involved bijection with the contratableaux associated
with certain Yamanouchi words,
which gives another combinatorial interpretation of the LR rule.
This combinatorial interpretation is in fact  different from
the one given by LR tableaux, which makes the matter even more
confusing.

 In a subsequent development, Knutson and Tao introduced~\cite{kt}
 the so called honeycombs, which are related to BZ triangles by
 a bijection that they sketch at the end.
The paper~\cite{gp} uses a related construction of ``web diagrams''
 for a different purpose.
The appendix in~\cite{kt} also introduces a
 different language of {\it hives}, which proved to be more flexible
 to restate the Knutson-Tao proof of saturation conjecture~\cite{bu}.

In the appendix to~\cite{bu}, Fulton described in a simple language
a bijection with a set of certain contratableaux, similar to that of
Carr\'e~\cite{car}.   As mentioned at the end of the appendix
(cf. also~\cite{ful}), the latter are in a well known bijection
with the classical LR tableaux.  Unfortunately, this bijection
is based on the Sch\" utzenberger involution, which is in fact
quite involved and goes beyond the scope of this paper.

\smallskip

Now, let us return to the linear maps establishing the bijections.
First, these maps show that the LR cones have the same
combinatorial structure.  Despite a visual difference
between definitions of LR tableaux, hives, and BZ triangles,
these combinatorial objects are essentially the same and should
be treated as equivalent.  In a sense, this varying nature of
these combinatorial interpretations of the LR coefficients makes
them ``more fundamental'' than others.

Let us mention here a ``local'' nature of the bijections
we present. A priori, linear maps $\vp: \rr^d \to \rr^d$ may
require $O(d^2)$ steps to perform.  In this case, however,
the local nature of bijections allows a~$O(d)$ computation.
This is especially striking when comparing the bijections
establishing the symmetries of the LR coefficients.
As observed previously, BZ triangles show nearly all the symmetries
of the LR coefficients, except for one: $c_{\mu, \nu}^\lambda =
c_{\nu, \mu}^\lambda$.  The latter seems to require $O(d^{3/2})$
operations and is significantly more difficult to perform.

Now, the idea of using integer points in
cones is a direct descendant of the earlier papers~\cite{gz,bz2}
and most recently has appeared in a context of integer
partitions~\cite{pak}.
While the fact that the linear maps between cones exist
at all may seem surprising, we do not claim to be the first
to establish that.  It is perhaps surprising that the
resulting linear maps are so simple and natural in this language.
We believe that this approach is perhaps more direct and
fruitful when compared to other more traditional combinatorial
techniques employed earlier (see above).

\smallskip

To conclude, let us describe  the structure of the paper.
We present in separate sections the LR tableaux, the hives
of Knutson and Tao, and the BZ triangles.  Along the way
we establish the bijections between these combinatorial
interpretations.
While the linear maps which produce these bijections are
easy to define, their proofs are not straightforward and
are delayed until the end of the paper.  We conclude with
the final remarks.

\vskip 2.5pc
{\large\bf 2\quad Littlewood-Richardson tableaux}

\vskip 1.5pc
Let $\lambda=\vector \lambda k$ be a partition of a positive integer $n$,
that is, a sequence of integers whose sum is $n$ and satisfy
$\lambda_1\ge \lambda_2\ge \cdots \ge \lambda_k \ge 0$.
Its diagram is the set of pairs of positive integers
$\{\, ( i, j) \mid 1 \le i \le k,\ 1 \le j \le \lambda_i \,\}$, which we also
denote by $\lambda$.
If $\mu$ is another partition and the diagram of $\mu$ is a
subset of the diagram of $\lambda$, in symbols $\mu \subseteq \lambda$, we
denote by $\lambda/ \mu$ the {\em skew diagram} consisting of the points in
$\lambda$ that are not in $\mu$, and by $\vert\lambda / \mu\vert$ its
cardinality.
It is customary to represent diagrams pictorially as a collection
of boxes \cite{ful, macd, stan}.
Any filling $T$ of a skew diagram $\lambda/ \mu$ with positive
integers, formally a map $T: \lambda/ \mu \flecha \natural$, will be called a
{\em Young tableau} or just a {\em tableau} of {\em shape} $\lambda/ \mu$.
A Young tableau $T$ is called {\em semistandard} if its rows are weakly
increasing from left to right and its columns are strictly increasing from top
to bottom.
The {\em \cont}of $T$ is the composition $\gamma (T)= \vector \gamma c$,
where $\gamma_i$ is the number of $i$'s in $T$.
The {\em word} of $T$, denoted by  $w(T)$ , is obtained
from $T$ by reading its entries from right to left, in successive rows,
starting with the top row and moving down.
For example, let
\[
D=\raisebox{.30em}{
\tableau{   &   &    & \ & \ & \ \\
                 &   & \ & \  \\
               \ & \ & \ & \  \\
               \ & \ & \ \\}}
\qquad {\rm and} \qquad
T=\raisebox{.30em}{
\tableau{   &    &    & 1 & 1 & 7 \\
                 &    & 1 & 4  \\
              1 & 4 & 5 & 7 \\
              3 & 5 & 7 \\}}
\]
then $D$ is a diagram of shape $( 6,4,4,3) / (3,2)$ and $T$ is a
semistandard tableaux of this shape, has content $(4,0,1,2,2,0,3)$ and
its word is $w(T) = 711417541753$.
Finally, a word $w = w_1\cdots w_k$ in the alphabet $1,\dots, n$ is called a
{\em lattice permutation} if for all $ 1\le  j \le k$  and all
$1\le i \le n -1$ the number of occurrences of $i$ in $w_1\cdots w_j$  is not
less than the number of occurrences of $i+1$ in $w_1\cdots w_j$.
A semistandard tableau $T$ of skew shape is called a
{\em Littlewood-Richardson} tableau if its word $w(T)$ is a lattice
permutation.
Note that the \cont of a \liri tableau is always a partition.
Given three partitions $\lambda$, $\mu$, $\nu$ such that
$\mu\subseteq\lambda$ and $\size\lambda=\size\mu+\size\nu$, we denote by
$c^\lambda_{\mu\,\nu}$ the number of \liri tableaux of shape $\lambda/\mu$ and
\cont $\nu$.
We will use the following example throughout the paper.
Let
\begin{equation}
\lambda=(23, 18, 15, 11, 8),\ \mu=(15, 9, 5, 2,0)\ {\rm and}\
\nu=(16,11,10,5,2),
\label{ejemplo}
\end{equation}
then the tableau in Figure \ref{fig:lrtab} is an example of a
\liri tableau of shape $\lambda/\mu$ and \cont $\nu$.
\begin{figure}[h]
\begin{center}
$
\tableau{   & & & & & & & & & & & & & & & 1 & 1 & 1 & 1 & 1 & 1 &
                                          1 & 1  \\
            & & & & & & & & & 1 & 1 & 1 & 1 & 2 & 2 & 2 & 2 & 2 \\
            & & & & & 1 & 1 & 2 & 2 & 2 & 2 & 2 & 3 & 3 & 3 \\
            & & 2 & 3 & 3 & 3 & 3 & 3 & 3 & 4 & 4 \\
            1 & 1 & 3 & 4 & 4 & 4 & 5 & 5 }$
\end{center}
\caption{\liri tableau}
\label{fig:lrtab}
\end{figure}

\vskip 2.5pc
{\large\bf 3\quad Littlewood-Richardson triangles}

\vskip 1.5pc
The {\em hive graph $\HG_k$ of size $k$} is a graph in the plane with
$\binom {k+2}{2}$ vertices arranged in a triangular grid consisting
of $k^2$ small equilateral triangles, as shown in
Figure~\ref{fig:hivegraph}.
\begin{figure}[h]
\begin{center}
\setlength{\unitlength}{1cm}
\begin{picture}(4,3.5)(0,.5)
\thicklines
\multiput(0,0)(1,0){5}{\circle*{.15}}
\multiput(.5,.866)(1,0){4}{\circle*{.15}}
\multiput(1,1.732)(1,0){3}{\circle*{.15}}
\multiput(1.5,2.598)(1,0){2}{\circle*{.15}}
\put(2,3.464){\circle*{.15}}
\drawline[1000](0,0)(2,3.464)(4,0)(0,0)
\drawline[1000](1,0)(2.5,2.598)
\drawline[1000](2,0)(3,1.732)
\drawline[1000](3,0)(3.5,.866)
\drawline[1000](1,0)(.5,.866)
\drawline[1000](2,0)(1,1.732)
\drawline[1000](3,0)(1.5,2.598)
\drawline(.5,.866)(3.5,.866)
\drawline(1,1.732)(3,1.732)
\drawline(1.5,2.598)(2.5,2.598)
\end{picture}
\end{center}
\caption{Hive graph $\HG_4$.}
\label{fig:hivegraph}
\end{figure}
Let $T_k$ denote the vector space of all labelings
$A=(a_{ij})_{0\le i\le j\le k}$ of the vertices of $\HG_k$ with real
numbers such that $a_{00}=0$.
We will write such labelings as triangular arrays of real
numbers in the way shown in Figure \ref{fig:triangarray}.
The dimension of $T_k$ is clearly $\binom{k+2}{2} -1$.
\begin{figure}[h]
\begin{center}
\setlength{\unitlength}{1cm}
\begin{picture}(3,2.5)(.3,.5)
\put(0,0){$\ent a03$}
\put(1,0){$\ent a13$}
\put(2,0){$\ent a23$}
\put(3,0){$\ent a33$}
\put(.5,.866){$\ent a02$}
\put(1.5,.866){$\ent a12$}
\put(2.5,.866){$\ent a22$}
\put(1,1.732){$\ent a01$}
\put(2,1.732){$\ent a11$}
\put(1.5,2.598){$\ent a00$}
\end{picture}
\end{center}
\caption{Triangular array of size 3.}
\label{fig:triangarray}
\end{figure}

We now proceed to explain how \liri tableaux can be coded in a
simple way as elements in $T_k$ satisfying certain inequalities.
A {\em \liri triangle of size $k$} is an element $A=(a_{ij})\in T_k$
that satisfies the following conditions:

\smallskip
(P) \ $a_{ij}\ge 0$, for all $1\le i < j\le k$.

(CS) \ $\sum_{p=0}^{i-1} a_{pj} \ge \sum_{p=0}^{i} a_{p\, j+1}$, for all
$1\le i\le j < k$.

(LR) \ $\sum_{q=i}^j a_{iq} \ge \sum_{q=i+1}^{j+1} a_{i+1\,q}$,
for all $1\le i\le j <k$.

\smallskip \noindent
Note that the inequality
\begin{equation}
\sum_{p=0}^{j} a_{pj} \ge \sum_{p=0}^{j+1} a_{p\, j+1},
\text{ for $1\le j < k$.}
\label{redundante}
\end{equation}
follows from (CS) with $i=j$  and (LR) with $i=j$;
also note that $a_{0j}$ and $a_{jj}$ could be negative.
We denote by $\lr k$ the cone of all \liri triangles in $T_k$,
and call it a {\em \liri cone}; this has the same dimension as $T_k$.
Also let $\pesos k$ denote the set of all $k$-tuples
$\lambda=\vector \lambda k$ of real numbers such that
$\lambda_1\ge \lambda_2\ge \cdots\ge \lambda_k$,
and $\size\lambda$ the sum of its entries, that is,
$\size \lambda = \sum_{i=1}^k \lambda_i$.
To each $A=(a_{ij}) \in \lr k$ we associate the following numbers:

\smallskip
(B1) \ $\mu_j = a_{0j}$, for all $1\le j\le k$.

(B2) \ $\lambda_j = \sum_{p=0}^j a_{pj}$, for all $1\le j\le k$.

(B3) \ $\nu_i = \sum_{q=i}^k a_{iq}$, for all $1\le i\le k$.

\smallskip \noindent
Then it follows from (P), (CS) and (LR) that the vectors
$\lambda=\vector\lambda k$, $\mu=\vector \mu k$ and
$\nu=\vector\nu k$ are in $\pesos k$ and that
$\size \lambda = \size \mu + \size \nu$.
We call $\tipo$ the {\em type} of $A$, and denote by  $\lrtipo k$
the set of all \liri triangles of type $\tipo$;
this is a convex polytope.
For example, let $\lambda$, $\mu$, $\nu$ be as in (\ref{ejemplo}), then the
triangle in Figure~\ref{fig:lrtri} is in $\lrtipo 5$.
\begin{figure}[h]
\begin{center}
\setlength{\unitlength}{1cm}
\begin{picture}(5,4.5)(0,.5)
\put(0,0){0}
\put(1,0){2}
\put(2,0){0}
\put(3,0){1}
\put(4,0){3}
\put(5,0){2}
\put(.5,.866){2}
\put(1.5,.866){0}
\put(2.5,.866){1}
\put(3.5,.866){6}
\put(4.5,.866){2}
\put(1,1.732){5}
\put(2,1.732){2}
\put(3,1.732){5}
\put(4,1.732){3}
\put(1.5,2.598){9}
\put(2.5,2.598){4}
\put(3.5,2.598){5}
\put(2,3.464){\makebox(0,0)[b]{\ 15}}
\put(3,3.464){8}
\put(2.5,4.33){0}
\end{picture}
\end{center}
\caption{\liri triangle of size 5}
\label{fig:lrtri}
\end{figure}

\vskip 1.5pc
Let $\lambda$, $\mu$, $\nu\in\pesos k$ be partitions, that is
$\lambda$, $\mu$ and $\nu$ have non-negative integer coefficients,
and suppose that $\size\lambda=\size\mu+\size\nu$.
To each \liri tableau $T$ of shape $\lambda/\mu$ and \cont $\nu$
we associate a triangular array $A_T=(a_{ij}) \in T_k$ by defining

\smallskip
(i) $a_{00}=0$, $a_{0j}=\mu_j$ for $1\le j\le k$, and

\smallskip
(ii) $a_{ij}$ equal to the number of $i$'s in row $j$ of $T$ for
$1\le i\le j\le k$.

\smallskip
\noindent
Note that the \liri triangle in Figure \ref{fig:lrtri} corresponds to the \liri
tableau in Figure \ref{fig:lrtab}.

\vskip 1.5pc
{\bf 3.1. Lemma.}
{\em Let $\lambda$, $\mu$, $\nu\in\pesos k$ be partitions such that
$\size\lambda=\size\mu+\size\nu$.
Then the correspondence $T \asocia A_T$ is a bijection between
the set of all \liri tableaux of shape $\lambda/\mu$ and \cont $\nu$
and the set of all \liri triangles of type $\tipo$ with integer
entries.
In particular $\lrtipo k$ has $\lrcoef$ integer points.}

\medskip
In effect, Lemma~3.1 translates combinatorics of \liri tableaux
into the language of integer points in polyhedra.  Various other
translations of this kind appear in the literature and are more
or less equivalent to ours.  A short ``verification style''
proof is given in Section~6.

\vskip 2.5pc
{\large\bf 4\quad Hives}

\begin{figure}[h]
\begin{center}
\setlength{\unitlength}{1cm}
\begin{picture}(8,1.6)
\put(0,0){\circle*{.15}}
\put(1,0){\circle*{.15}}
\put(.5,.866){\circle*{.15}}
\put(1.5,.866){\circle*{.15}}
\put(3.5,.433){\circle*{.15}}
\put(4.5,.433){\circle*{.15}}
\put(4,-.433){\circle*{.15}}
\put(4,1.299){\circle*{.15}}
\put(6.5,.866){\circle*{.15}}
\put(7.5,.866){\circle*{.15}}
\put(7,0){\circle*{.15}}
\put(8,0){\circle*{.15}}
\drawline[1000](0,0)(1,0)(1.5,.866)(.5,.866)(0,0)
\drawline[1000](3.5,.433)(4,1.299)(4.5,.433)(4,-.433)(3.5,.433)
\drawline[1000](6.5,.866)(7.5,.866)(8,0)(7,0)(6.5,.866)
\end{picture}
\end{center}
\caption{Types of rhombi in a hive graph.}
\label{fig:rhombi}
\end{figure}
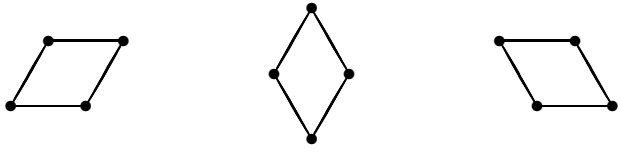

\vskip 1.5pc
The hive graph $\HG_k$ of size $k$ is divided into $k^2$ small
equilateral triangles.
Each two adjacent such triangles form a rhombus with two obtuse
angles and two acute angles.
There are three types of rhombi: tilted to the right, vertical and
tilted to the left.
They are shown in Figure \ref{fig:rhombi}.

A {\em hive of size $k$} is a labeling $H=\arreglo h$ of the vertices of
the hive graph $\HG_k$ with real numbers such that for each rhombus the
sum of the labels at obtuse vertices is bigger than or equal to the sum
of the labels at acute vertices; equivalently, $H=(h_{ij})$ satisfies the
following inequalities:

\smallskip
(R) \ $h_{i\, j} - h_{i\, j-1}\ge h_{i-1\, j} - h_{i-1\,j-1}$,
for $1\le i < j \le k$.

(V) \ $h_{i-1\,j} - h_{i-1\, j-1} \ge h_{i\, j+1}-h_{i\, j}$,
for $1\le i \le j < k$.

(L) \ $h_{i\,j} - h_{i-1\, j} \ge h_{i+1\, j+1}-h_{i\, j+1}$,
for $1\le i \le j <k$.

\smallskip
We denote by $\hive k$ the cone of all hives of size $k$ that satisfy
the extra condition $h_{00}=0$, and call it a {\em hive cone}.
As we did for \liri triangles, we associate to each hive
$H=(h_{ij})\in \hive k$ numbers:

\smallskip
(B1$^\prime$) \ $\mu_j = h_{0\, j} - h_{0\, j-1}$, for $1\le j\le k$.

(B2$^\prime$) \ $\lambda_j = h_{j\, j} - h_{j-1\, j-1}$, for $1\le j\le k$.

(B3$^\prime$) \ $\nu_i = h_{i\, k}- h_{i-1\, k}$, for $1\le i\le k$.

\smallskip
Then it follows from (R), (V) and (L) that the vectors
$\lambda=\vector\lambda k$, $\mu=\vector \mu k$ and
$\nu=\vector\nu k$ are in $\pesos k$ and that
$\size \lambda = \size \mu + \size \nu$.
For example,
\[
\mu_j = \ent h0j - \ent h0{j-1} \ge \ent h1{j+1} - \ent h1j \ge
\ent h0{j+1} - \ent h0j = \mu_{j+1}.
\]
We call $\tipo$ the {\em type} of $A$, and denote by  $\hivetipo k$
the set of all hives of type $\tipo$; this is a convex polytope.
For example, let $\lambda$, $\mu$ and $\nu$ be as in (\ref{ejemplo}),
then the triangle in Figure~\ref{fig:hivetri} is in $\hivetipo 5$.
\begin{figure}[h]
\begin{center}
\setlength{\unitlength}{1cm}
\begin{picture}(5,4.5)(0,.5)
\put(0,0){31}
\put(1,0){47}
\put(2,0){58}
\put(3,0){68}
\put(4,0){73}
\put(5,0){75}
\put(.5,.866){31}
\put(1.5,.866){45}
\put(2.5,.866){56}
\put(3.5,.866){65}
\put(4.5,.866){67}
\put(1,1.732){29}
\put(2,1.732){43}
\put(3,1.732){53}
\put(4,1.732){56}
\put(1.5,2.598){24}
\put(2.5,2.598){36}
\put(3.5,2.598){41}
\put(2,3.464){15}
\put(3,3.464){23}
\put(2.5,4.33){\ 0}
\end{picture}
\end{center}
\caption{Hive of size 5}
\label{fig:hivetri}
\end{figure}

For any positive integer $k$, we define a linear map
\, $\Phi_k \colon T_k \flecha T_k$ \,
by
$$\Phi_k(\ent aij)=(\ent hij), \ \ \text{where}  \ \
\ent hij =\sum_{p=0}^i\sum_{q=p}^j  \ent apq \,.$$
Note that the hive in Figure \ref{fig:hivetri} is the image under
$\Phi_5$ of the \liri triangle in Figure \ref{fig:lrtri}.
We have the following theorem.

\vskip 1.5pc
{\bf 4.1. Theorem.}
{\em The map $\Phi_k$ defined above
is a volume preserving linear operator
which maps $\lr k$ bijectively onto $\hive k$, and
$\lrtipo k$ onto $\hivetipo k$, for all $\lambda$,
$\mu$, $\nu\in\pesos k$.}

\bigskip
As mentioned in the introduction, the proof can be found
in section~6.  Let us mention here two important corollaries.
For any polytope~$P$ let $e(P)$ denote the number of integer
points in $P$.

\vskip 1pc
{\bf 4.2. Corollary.}
{\em $e(\hivetipo k)=c^\lambda_{\mu\,\nu}$, for all $\lambda$,
$\mu$, $\nu\in\pesos k$ with non-negative integer coefficients.}

\vskip 1pc
{\bf 4.3. Corollary.}
{\em ${\rm Vol}(\hivetipo k)= {\rm Vol}(\lrtipo k)$, for all $\lambda$,
$\mu$, $\nu\in\pesos k$.}


\vskip 2.5pc
{\large\bf 5\quad \beze triangles}

\vskip 1.5pc
For any integer $k\ge 1$ we construct a graph $\BZ_k$ from the hive graph
$\HG_{k+1}$ in the following way:
Its vertices are the middle points of the edges of the hive graph that
do not lie on the boundary, and their edges are those joining pairs of
middle points on edges lying on small triangles of $\HG_{k+1}$,
see Figure~\ref{fig:bz-graph}.
We call $\BZ_k$ the \beze graph of size $k$.
\begin{figure}[h]
\begin{center}
\setlength{\unitlength}{.7cm}
\begin{picture}(6.5,6.5)(-1.5,-.366)
\thicklines
\multiput(0,0)(1,0){6}{\circle*{.15}}
\multiput(.5,.866)(2,0){3}{\circle*{.15}}
\multiput(1,1.732)(1,0){4}{\circle*{.15}}
\multiput(1.5,2.598)(2,0){2}{\circle*{.15}}
\multiput(2,3.464)(1,0){2}{\circle*{.15}}
\put(2.5,4.33){\circle*{.15}}
\drawline[1000](-1.5,-.866)(6.5,-.866)(2.5,6.062)(-1.5,-.866)
\drawline[1000](-.5,.866)(5.5,.866)(4.5,-.866)(1.5,4.33)(3.5,4.33)
(.5,-.866)(-.5,.866)
\drawline[1000](.5,2.598)(4.5,2.598)(2.5,-.866)(.5,2.598)
\end{picture}
\qquad\qquad\qquad
\begin{picture}(6,5.5)(-1,-.366)
\thicklines
\multiput(0,0)(1,0){6}{\circle*{.15}}
\multiput(.5,.866)(2,0){3}{\circle*{.15}}
\multiput(1,1.732)(1,0){4}{\circle*{.15}}
\multiput(1.5,2.598)(2,0){2}{\circle*{.15}}
\multiput(2,3.464)(1,0){2}{\circle*{.15}}
\put(2.5,4.33){\circle*{.15}}
\drawline[1000](0,0)(2.5,4.33)(5,0)(0,0)
\drawline[1000](1,0)(.5,.866)
\drawline[1000](2,0)(3.5,2.598)
\drawline[1000](3,0)(1.5,2.598)
\drawline[1000](4,0)(4.5,.866)
\drawline(1,1.732)(4,1.732)
\drawline(2,3.464)(3,3.464)
\end{picture}
\end{center}
\caption{Hive graph $\HG_4$ and the corresponding \beze graph $\BZ_3$.}
\label{fig:bz-graph}
\end{figure}
The vertices of the \beze graph are partitioned into disjoint blocks of
cardinality three, each block corresponding to a small equilateral triangle;
these triangles are distributed in the graph: one on the first (top) level,
two on the second level, three on the third level, and so on.
Let $V_k$ denote the vector space of all labelings
$X = (\ent xij, \ent yij, \ent zij)_{1\le i\le j\le k}$ of
$\BZ_k$ with real numbers.
The labelings are carried out in such a way that the vertices of the
$i$-th triangle on the $j$-th level are labeled with $\ent xij$,
$\ent yij$, $\ent zij$ as indicated in Figure~\ref{fig:bz-label}.
\begin{figure}[h]
\begin{center}
\setlength{\unitlength}{1cm}
\begin{picture}(6,4.33)(-.2,.6)
\put(0,0){$\ent y13$}
\put(1,0){$\ent z13$}
\put(2,0){$\ent y23$}
\put(3,0){$\ent z23$}
\put(4,0){$\ent y33$}
\put(5,0){$\ent z33$}
\put(.5,.866){$\ent x13$}
\put(2.5,.866){$\ent x23$}
\put(4.5,.866){$\ent x33$}
\put(1,1.732){$\ent y12$}
\put(2,1.732){$\ent z12$}
\put(3,1.732){$\ent y22$}
\put(4,1.732){$\ent z22$}
\put(1.5,2.598){$\ent x12$}
\put(3.5,2.598){$\ent x22$}
\put(2,3.464){$\ent y11$}
\put(3,3.464){$\ent z11$}
\put(2.5,4.33){$\ent x11$}
\end{picture}
\end{center}
\caption{Labeling of $\BZ_3$}
\label{fig:bz-label}
\end{figure}
The dimension of $V_k$ is $3\binom{k+1}{2}$.
Note that the labels $\ent yij$, $\ent zij$, $\ent x{i+1}{j+1}$,
$\ent y{i+1}{j+1}$, $\ent zi{j+1}$, $\ent xi{j+1}$ form an hexagon
for each $1\le i \le j <k$ and hence there are $\binom{k}{2}$ hexagons
in $\BZ_k$.
We will be interested in the subspace $W_k$ of $V_k$ consisting of all
labelings such that for each hexagon in $\BZ_k$ the sum of the labels in
each edge equals the sum of the labels of the diametrically opposite edge,
that is

\smallskip
(BZ1) \ $\ent yij + \ent zij = \ent y{i+1}{j+1} + \ent zi{j+1}$,
for all $1\le i\le j < k$.

(BZ2) \ $\ent xi{j+1} + \ent yij = \ent x{i+1}{j+1} + \ent y{i+1}{j+1}$,
for all $1\le i\le j < k$.

(BZ3) \ $\ent xi{j+1} + \ent zi{j+1} = \ent x{i+1}{j+1} + \ent zij$,
for all $1\le i\le j < k$.

\smallskip
Observe that any of these three equalities follows from the other two.

\vskip 1.5pc
{\bf 5.1. Lemma.}
{\em The vector space $W_k$ has dimension $\frac{1}{2} k (k+5)=
{\rm dim}\, T_{k+1}-2$.}

\bigskip
A {\em \beze triangle of size k} is any labeling of $\BZ_k$ in $W_k$
with non-negative entries.
Let $\bz k$ denote the cone of all \beze triangles of size $k$.
Let $\lambda$, $\mu$, $\nu\in\pesos {k+1}$, then we say that a \beze triangle
is of type $\tipo$ is it satisfies the following conditions:

(B1$^{\prime\prime}$) \ $\ent x1j + \ent y1j =
                      \mu_j - \mu_{j+1}$, for $1\le j\le k$.

(B2$^{\prime\prime}$) \ $\ent xjj + \ent zjj =
                       \lambda_j - \lambda_{j+1}$, for $1\le j\le k$.

(B3$^{\prime\prime}$) \ $\ent yik + \ent zik =
                        \nu_i - \nu_{i+1}$, for $1\le i\le k$.

\smallskip
Note that, in contrast to \liri triangles and hives, a \beze triangle
has many different types.
Let $\bztipo k$ denote the set of all \beze triangles of type $\tipo$;
this is a convex polytope.
For example, let $\lambda$, $\mu$ and $\nu$ be as in (\ref{ejemplo}),
then the triangle in Figure \ref{fig:bztri} is in $\bztipo 4$.
Here the $\ent xij$'s are written with roman numerals, the $\ent yij$'s
by {\bf boldface} numerals and the $\ent zij$'s by {\it italic\/} numerals.
\begin{figure}[h]
\begin{center}
\setlength{\unitlength}{1cm}
\begin{picture}(7,6)(.1,.4)
\put(0,0){\bf 2}
\put(1,0){\it 3}
\put(2,0){\bf 0}
\put(3,0){\it 1}
\put(4,0){\bf 1}
\put(5,0){\it 4}
\put(6,0){\bf 3}
\put(7,0){\it 0}
\put(.5,.866){0}
\put(2.5,.866){0}
\put(4.5,.866){0}
\put(6.5,.866){3}
\put(1,1.732){\bf 0}
\put(2,1.732){\it 3}
\put(3,1.732){\bf 1}
\put(4,1.732){\it 1}
\put(5,1.732){\bf 6}
\put(6,1.732){\it 1}
\put(1.5,2.598){3}
\put(3.5,2.598){4}
\put(5.5,2.598){3}
\put(2,3.464){\bf 2}
\put(3,3.464){\it 2}
\put(4,3.464){\bf 5}
\put(5,3.464){\it 2}
\put(2.5,4.33){2}
\put(4.5,4.33){1}
\put(3,5.196){\bf 4}
\put(4,5.196){\it 3}
\put(3.5,6.062){2}
\end{picture}
\end{center}
\caption{\beze triangle of size 4.}
\label{fig:bztri}
\end{figure}

\vskip 1.5pc
For any integer $k \ge 2$, we define a linear map
\, $\Psi_k \colon T_k \flecha W_{k-1}$ \,
by setting \,
$\Psi_k(\ent hij) = (\ent xij, \ent yij, \ent zij )$ \,
where
\begin{eqnarray*}
\ent xij & = & \ent hij + \ent h{i-1}j - \ent h{i-1}{j-1} -\ent hi{j+1}, \\
\ent yij & = & \ent h{i-1}j + \ent hi{j+1} - \ent hij - \ent h{i-1}{j+1}, \\
\ent zij & = & \ent hij + \ent hi{j+1} - \ent h{i-1}j - \ent h{i+1}{j+1},
\end{eqnarray*}
for all $1\le i\le j < k$.
Note that the values of the $\ent xij$'s, $\ent yij$'s and $\ent zij$'s
are obtained by taking, respectively, the differences of the inequalities
(V), (R) and (L) used to define hives.
It should be remarked that the $\ent yij$'s are obtained from (R) by
adding one to $j$.
It is straightforward to check that the image of $\Phi_k$ is contained
in $W_{k-1}$.
The composition \,
$\Psi_k \circ \Phi_k \colon T_k \flecha W_{k-1} $ \,
has also a nice description:\, $\Psi_k \circ \Phi_k(\ent aij)=
(\ent xij, \ent yij, \ent zij )$\, with
\begin{eqnarray}
\ent xij & = & \textstyle{\sum_{p=0}^{i-1} \ent apj -
\sum_{p=0}^i \ent ap{j+1} }, \notag\\
\ent yij & = & \ent ai{j+1}, \label{imagen}\\
\ent zij & = & \textstyle{\sum_{q=i}^j \ent aiq -
\sum_{q=i+1}^{j+1} \ent a{i+1}q}, \notag
\end{eqnarray}
for all $1\le i\le j < k$.
Again, the values of the $\ent xij$'s, $\ent yij$'s and $\ent zij$'s
are obtained by taking, respectively, the differences of the left and
right hand sides in the inequalities (CS), (P) and (LR) used to
define \liri triangles.
For example, the \beze triangle in Figure~\ref{fig:bztri} is
the image under $\Psi_5$ of the hive in Figure~\ref{fig:hivetri} and
the image under $\Psi_5\circ\Phi_5$ of the \liri triangle in
Figure~\ref{fig:lrtri}.
Note that the boldface numerals in Figure~\ref{fig:bztri} are contained
in the \liri triangle from Figure~\ref{fig:lrtri}.

\vskip 1.5pc
{\bf 5.2. Theorem.}
{\em The linear operator $\Psi_k\circ\Phi_k$ maps $\lr k$ surjectively
onto $\bz {k-1}$, and $\lrtipo k$ bijectively onto $\bztipo {k-1}$, for any
$\lambda$, $\mu$, $\nu\in \pesos k$.}

\vskip 1.5pc
{\bf 5.3. Corollary.}
{\em The linear operator $\Psi_k$ maps $\hive k$ surjectively onto $\bz {k-1}$,
and $\hivetipo k$ bijectively onto $\bztipo {k-1}$, for any
$\lambda$, $\mu$, $\nu\in \pesos k$.}

\vskip 1.5pc
{\bf 5.4. Corollary.}
{\em $e(\bztipo {k-1})= \lrcoef$, for any $\lambda$, $\mu$,
$\nu\in \pesos k$ with non-negative integer coefficients.}

\vskip 1.5pc
It will follow from Lemma 6.1 and the proof of Theorem~5.2 that the cones
$\lr k$ and $\hive k$ are isomorphic to $\bz{k-1} \times \real^2$.
One can embed the cone $\bz{k-1}$ into $\lr k$ in the following way:
For any $k\ge 2$, let \, $\Omega _k \colon W_{k-1} \flecha T_k$ \,
be the linear operator defined by
$\Omega_k(\ent xij, \ent yij, \ent zij )= (\ent aij)$ where
\begin{eqnarray}
\ent a0j & = &  \textstyle{\sum_{l=j}^{k-1} \ent x1l + \ent y1l}
\ \text{ for } 1\le j < k, \ \text{ and }\ \ent a0k = 0, \notag\\
\ent aij & = & \ent yi{j-1}, \ \text{ for } 1\le i < j\le  k,
\label{embedding} \\
\ent ajj & = &  \textstyle{\sum_{l=j}^{k-1} \ent zll}\ \text{ for }
1\le j < k,  \ \text{ and }\ \ent akk = 0. \notag
\end{eqnarray}
Then we have:

\bigskip
{\bf 5.5. Theorem.}
{\em The linear operator $\Omega_k$ defined above maps
$\bz {k-1}$ injectively into $\lr k$, and $\bztipo {k-1}$ bijectively
onto $\lrtipo k$ for any $\lambda$, $\mu$, $\nu\in\pesos k$ such that
$\mu_k=0$ and $\nu_k=0$.}

\vskip 2.5pc
{\large\bf 6 \quad Proof of results}

\vskip 1.5pc

\begin{proof}[\bf \indent Proof  of Lemma 3.1]
Let $T$ be a \liri tableau of shape $\lambda/\mu$ and \cont $\nu$,
then $A_T$ satisfies (P) by definition.
Since $T$ has strictly increasing columns (CS) follows,
and since $w(T)$ is a lattice permutation, $A_T$ satisfies (LR).
It is also clear that $A_T$ is of type $\tipo$.
Conversely, for any \liri triangle $A=(a_{ij})$ in $\lrtipo k$
with integer entries, we define a tableau $T_A$ of shape $\lambda/\mu$
by placing in row $j$, in weakly increasing order, $a_{ij}$ $i$'s for
each $i$ and $j$.
It is routine to check that $T$ is a \liri tableau of shape
$\lambda/\mu$ and \cont $\nu$, and that both constructions are
inverses of each other.
Here we use that (\ref{redundante}) follows from (CS) and (LR).
\end{proof}

\vskip 1pc
\begin{proof}[\bf \indent Proof of Theorem 4.1]
Let $\{E_{ij}\}$ be the canonical basis of $T_k$,
that is $E_{ij}=\left(e^{ij}_{pq}\right)$, where
\begin{equation*}
e^{ij}_{pq}=
\begin{cases}
1, & \text{if $p=i$ and $q=j$;} \\
0, & \text{otherwise.}
\end{cases}
\end{equation*}
We order it according to the lexicographic order of the subindices,
that is,
\[
{\cal B}=\{ E_{01},E_{02}, \dots, E_{0k}, E_{11}, \dots, E_{1k},
\dots, E_{kk}\}.
\]
The matrix of $\Phi_k$ with respect to ${\cal B}$ is lower
triangular with ones on the main diagonal, therefore it has determinant
one, is volume preserving, and maps $\entero^{\binom{k+2}{2}-1}$
bijectively onto $\entero^{\binom{k+2}{2}-1}$.
The inverse of $\Phi_k$ is given by $\Phi_k^{-1}(\ent hij) = (\ent aij)$
where
\[
\ent aij =
\begin{cases}
\ent h0j - \ent h0{j-1}, & \text{if $i=0$ and $1\le j\le k$.} \\
\ent hjj - \ent h{j-1}j, & \text{if $1\le i=j\le k$.} \\
\ent hij - \ent hi{j-1} - \ent h{i-1}j + \ent h{i-1}{j-1},
& \text{if $1\le i<j\le k$.}
\end{cases}
\]
Let $(\ent aij)\in\lr k$ and $(\ent hij) = \Phi_k(\ent aij)$,
then we have
\[
\ent hst - \ent hs{t-1} = \sum_{p=0}^s \ent apt \ \text{ and } \
\ent h{s+1}t - \ent hst = \sum_{q=s+1}^t \ent a{s+1}q,
\]
for $0\le s<t\le k$.
It is straightforward, using these two identities,
to check that $(\ent aij)$ satisfies (P), (CS) or (LR), respectively,
if and only if $(\ent hij)$ satisfies (R), (V) or (L), respectively;
therefore $\Phi_k(\lr k) = \hive k$.
Also, it is straightforward to check that $(\ent aij)$ and
$(\ent hij)$ have the same type; therefore
$\Phi_k(\lrtipo k) = \hivetipo k$, for all $\lambda$, $\mu$
and $\nu\in \pesos k$.
\end{proof}

\vskip1.5pc
\begin{proof}[\bf \indent Proof of Lemma 5.1]
We form a system of linear equations by taking, for each $1\le i\le j <k$,
that is, for each hexagon in $\BZ_k$, equations (BZ2) and (BZ3).
Then, after arranging the variables in the order $\ent x11$, $\ent y11$,
$\ent z11$, $\ent x12$, $\ent y12$, $\ent z12$, $\ent x22, \dots, \ent zkk$,
we easily check that the matrix of coefficients of the system is in
echelon form and has rank $2\binom{k}{2}$.
Thus ${\rm dim}\, W_k=3\binom{k+1}{2} -2 \binom{k}{2}=\frac{1}{2} k(k+5)$.
\end{proof}

\vskip 1.5pc
Before we prove Theorem 5.2, let us prove the following lemma.

\bigskip
{\bf 6.1. Lemma.}
{\em The linear operators $\Psi _k$ and $\Psi_k\circ\Phi_k$ are surjective.
Moreover, equations~(\ref{preimagen}) give a full description of
$(\Psi_k\circ\Phi_k)^{-1}(X)$ for any $X\in W_{k-1}$.}
\begin{proof}[\bf \indent Proof]
It is enough to show that $\Psi_k\circ\Phi_k$ is surjective.
Let $X=(\ent xij, \ent yij, \ent zij)\in W_{k-1}$.
For each $s$, $t\in\real$ we define an element
$\ent Ast=(\ent aij)\in T_k$ by
\begin{eqnarray}
\ent a0j & = & s + \textstyle{\sum_{l=j}^{k-1} \ent x1l + \ent y1l}
\ \text{ for } 1\le j < k, \ \text{ and }\ \ent a0k = s, \notag\\
\ent aij & = & \ent yi{j-1}, \ \text{ for } 1\le i < j\le  k,
\label{preimagen} \\
\ent ajj & = & t + \textstyle{\sum_{l=j}^{k-1} \ent zll}\ \text{ for }
1\le j < k,  \ \text{ and }\ \ent akk = t. \notag
\end{eqnarray}
Then it follows from repeated application of (BZ2) and (BZ1) that
$\Psi_k\circ\Phi_k(\ent Ast)=X$.
The last statement follows from the identity
${\rm dim}\, T_k= {\rm dim}\, W_{k-1} + 2$.
\end{proof}

\vskip 1.5pc
\begin{proof}[\bf \indent Proof of Theorem 5.2]
It follows from (\ref{imagen}) that $\Psi_k\circ\Phi_k(\lr k) = \bz {k-1}$;
and it follows from (\ref{imagen}) and (B1)-(B3) that $A$ and
$\Psi_k\circ\Phi_k(A)$ have the same type, for any $A\in\lr k$, thus
$\Psi_k\circ\Phi_k(\lrtipo k)=\bztipo k$.
The last claim follows from the remark that different elements
in the preimage of an $X\in\bztipo k$ have different types.
\end{proof}

\vskip 1.5pc
\begin{proof}[\bf \indent Proof of Theorem 5.5.]
It follows from (\ref{imagen}), (\ref{embedding}) and the proof of
Lemma 6.1 that $\Psi_k \circ \Phi_k \circ \Omega_k$ is the identity
map on $W_{k-1}$, and that $\Omega(\bz{k-1})\subseteq \lr k$.
The last statement follows from Theorem~5.2.
\end{proof}

\vskip 2.5pc
{\large\bf 7 \quad Final remarks}

\vskip 1.5pc

Let us start with the complexity issues.  Recall that the LR triangles,
hives, and BZ triangles, all of size~$k$, are given by $\theta(k^2)$
entries.  As defined, maps $\Phi^{-1}$ and
$\Psi$ require only a constant
number of operations per entry, i.e. have $O(k^2)$ complexity.  It is
an easy exercise in dynamic programming to show that $\Phi$ and
$\Psi^{-1}$ have the same complexity, linear in the input.

The complexity $O(k^2)$ is in stark contrast with the $O(k^3)$
complexity required by the jeu-de-taquin and Sch\" utzenberger
involution (cf.~\cite{ful,stan}).  This explains why Fulton's
map in~\cite{bu} has the same complexity.  In fact, Fulton reworks
the bijection of Carr\'e~\cite{car} which establishes a combinatorial
map $\Upsilon: e(\lrtipo k) \to e(\carre k)$.  As we mentioned in the
introduction and will reiterate below, there is no linear map
establishing the symmetry $\carre k \to \hivetipo k$.  One can use
a more complicated map called tableaux switching to demonstrate
this symmetry~\cite{bs} (see also~\cite{lee}).

\smallskip

Now, the symmetries of the LR coefficients are quite intriguing in
a sense that all but one of them can be established by simple means.
If one operates with LR tableaux, one simply has to map them into
BZ triangles (which takes $O(k^2)$ steps), perform the symmetry,
and return back to LR tableaux (which takes $O(k^2)$ steps again).
For the remaining $\mu \to \nu$ symmetry several authors found
an explicit map (in different languages)~\cite{aze1,aze2,lee,tao} but
all of them use $O(k^3)$ steps.
It would be interesting to prove the lower bound
$\Omega(k^3)$ but we are doubtful such result is feasible
at the moment.  What one can show, however, is that this `last'
symmetry cannot be performed by a linear map already for~$k=4$.
We leave this statement as an interesting exercise to the reader,
in the hope that further results will be found in this direction.

\newpage

\vskip 2.5pc
{\bf Acknowledgements}

\smallskip
\noindent
We are grateful to Oleg Gleizer, Michael Kleber,
Alex Postnikov and Terry Tao
for interesting conversations and helpful remarks.
We thank Olga Azenhas and
Christophe Carr\'e for the help with the literature.

The first author was supported by NSA and NSF.
The second author is grateful to Richard Stanley for
his support in organization of the sabbatical visit
to~MIT.


\vskip1.cm
 \smallskip
 \noindent
 {\bf Keywords}:  Young tableaux, Littlewood-Richardson rule,
 Berenstein-Zelevinsky triangles, Knutson-Tao hives


\vskip 2.5pc
{\small
\begin{thebibliography}{99}

\bibitem{aze1} O.~Azenhas, Littlewood-Richardson fillings
and their symmetries,
{\em Textos de Matem\'atica} S\'erie B, {\bf 19} (1999), 81--92

\bibitem{aze2} O.~Azenhas, On an involution on the set of
Littewood-Richardson tableaux and the hidden commutativity,
Pr\'e-publica\c{c}\~{o}es do Departamento de Matem\'atica
da Uni\-versidade de Coimbra 00-27 (2000); available from
{\tt http://dingo.mat.uc.pt/}$\widetilde\, $ {\tt cmuc/publicline.php?lid=1}

\bibitem{bs} G.~Benkart, F.~Sottile, J.~Stroomer,
Tableau switching: algorithms and applications
{\em J. Combin. Theory}, Ser. A  \, {\bf 76} (1996), 11--43.

\bibitem{bz1}  A.~D. Berenstein and A.~V. Zelevinsky,
Involutions on Gelfand-Tsetlin schemes and
multiplicities in skew ${\rm GL}\sb n$-modules (in Russian)
{\em Dokl. Akad. Nauk SSSR } {\bf 300}  (1988),  1291--1294

\bibitem{bz2}  A.~D. Berenstein and A.~V. Zelevinsky,
Tensor product multiplicities and convex polytopes in partition space
{\em J. Geom. Phys.} {\bf 5}  (1988),  no. 3, 453--472

\bibitem{beze} A.~D. Berenstein and A. Zelevinsky,
Triple multiplicities for $sl(r+1)$ and the spectrum of the exterior
algebra of the adjoint representation,
{\em  J. Algebraic Combin.} {\bf 1} (1992), 7-22.

\bibitem{bu} A. Buch,
The saturation conjecture (after A. Knutson and T. Tao),
{\em Enseign. Math.} {\bf 46} (2000), 43-60.

\bibitem{car} C. Carr\'e,
The rule of \liri in a construction of Berenstein-Zelevinsky,
{\em Internat. J. Algebra Comput.} {\bf 1} (1991), 473-491.

\bibitem{ful} W. Fulton,
{\em Young Tableaux},
London Math. Soc. Student Texts 35,
Cambridge Univ. Press 1997.

\bibitem{gz}
I.~M. Gelfand and A.~V. Zelevinsky,
Multiplicities and regular bases for ${\rm gl}\sb n$ (in Russian),
{\em Group-theoretic methods in physics}, Vol. 2
(J\=urmala, 1985), 22--31, "Nauka", Moscow, 1986.

\bibitem{gp} O. Gleizer, A. Postnikov,
Littlewood-Richardson coefficients via Yang-Baxter equation.
{\em Internat. Math. Res. Notices} (2000), no. 14, 741--774

\bibitem{kb}
A.~N. Kirillov and A.~D. Berenstein,
Groups generated by involutions, Gelfand-Tsetlin patterns,
and combinatorics of Young tableaux,
{\em Algebra i Analiz} {\bf 7} (1995), 92--152

\bibitem{kt} A. Knutson and T. Tao,
The honeycomb model of GL$_n(\complejo)$ tensor products I:
Proof of the saturation conjecture,
{\em J. Amer. Math. Soc.} {\bf 12} (1999), 1055-1090.

\bibitem{lee} M.~A.~A. van Leeuwen,
The Littlewood-Richardson rule, and related combinatorics, in
{\em Interaction of combinatorics and representation theory},
95--145,  {\em MSJ Mem.}, {\bf 11}, Math. Soc. Japan, Tokyo, 2001

\bibitem{macd} I.~G. Macdonald,
{\em Symmetric functions and Hall polynomials}, 2nd. ed.,
Oxford University Press, 1995.

\bibitem{pak} I. Pak,
Partition Identities and Geometric Bijections, {\em Proc. A.M.S.},
to appear (2002)

\bibitem{stan} R.P. Stanley,
{\em Enumerative Combinatorics}, vol. 2,
Cambridge Studies in Advanced Mathematics 62. Cambridge Univ. Press, 1999.

\bibitem{tao} T. Tao, personal communication

\bibitem{zel}
A. Zelevinsky, Littlewood-Richardson semigroups, in
{\em New perspectives in algebraic combinatorics} (Berkeley, CA, 1996--97),
337--345, Math. Sci. Res. Inst. Publ., {\bf 38},
Cambridge Univ. Press, Cambridge, 1999.

\end{thebibliography}
}
\end{document}